\documentclass[11pt]{amsart}

\usepackage{amssymb}
\usepackage{epsfig}
\usepackage{a4wide}

\newtheorem{theorem}{Theorem}

\newtheorem{lemma}[theorem]{Lemma}
\newtheorem{corollary}[theorem]{Corollary}

\newcommand{\ts}{\hspace{0.5pt}}
\newcommand{\CC}{\mathbb{C}\ts}
\newcommand{\RR}{\mathbb{R}\ts}
\newcommand{\QQ}{\mathbb{Q}\ts}
\newcommand{\FF}{\mathbb{F}\ts}

\newcommand{\cE}{\mathcal{E}}

\newcommand{\cR}{\mathcal{R}}

\def\mod#1{\,({\rm mod\ }#1)}
\def\proof{\noindent{\sc Proof.}\hskip 5pt}
\def\endproof{\hfill\vbox{\hrule
    \hbox{\vrule\kern4pt\vbox{\kern4pt
    \kern4pt}\kern4pt\vrule}\hrule}\bigskip}

\begin{document}

\title[Combinatorial model for reversible maps]
{A combinatorial model for reversible \\rational maps over
finite fields}

\author{John A.~G.~Roberts}
\address{School of Mathematics and Statistics,
University of New South Wales,
Sydney, NSW 2052, Australia}
\email{jag.roberts@unsw.edu.au}
\urladdr{http://www.maths.unsw.edu.au/\~{}jagr}

\author{Franco Vivaldi}
\address{School of Mathematical Sciences, Queen Mary,
University of London,
London E1 4NS, UK}
\email{f.vivaldi@maths.qmul.ac.uk}
\urladdr{http://www.maths.qmul.ac.uk/\~{}fv}

\begin{abstract}
We study time-reversal symmetry in dynamical systems with finite phase 
space, with applications to birational maps reduced over finite fields. 
For a polynomial automorphism with a single family of reversing symmetries, a 
universal (i.e., map-independent) distribution function $\cR(x)=1-e^{-x}(1+x)$
has been conjectured to exist, for the normalized cycle lengths of the 
reduced map in the large field limit \cite{RobertsVivaldi:05}.
We show that these statistics correspond to those of a composition of 
two random involutions, having an appropriate number of fixed points. 
This model also explains the experimental observation that, 
asymptotically, almost all cycles are symmetrical, and that the probability 
of occurrence of repeated periods is governed by a Poisson law.
\end{abstract}

\maketitle

\section{Introduction}

This paper is concerned with time-reversal symmetry in dynamical
systems with finite phase space, with applications to rational 
maps over finite fields. The concept of reversibility originated in
the theory of smooth maps and flows on manifolds \cite{DeVog,Finn,Devaney}. 
In this paper, a map $L$ is said to be {\it reversible\/} if it is 
the composition of two involutions $G$ and $H=L\circ G$.
The involution $G$ conjugates the map to its inverse, namely
\begin{equation} \label{eq:RevSym}
         G \circ L \circ G^{-1} \; = \;  L^{-1}\,.
\end{equation}
(For background information on reversibility and its generalization, see
\cite{LambRoberts,MacKay,RobertsQuispel}, and references therein.)
The nature of the fixed sets of $G$ and $H$, respectively $Fix\,G$ 
and $Fix\,H$, plays a key dynamical role in organizing the dynamics of
a reversible map.

The reversibility property can be interpreted in purely algebraic terms, 
and in particular it applies to reversible algebraic maps $L$ of any dimension 
$n$ over any field $K$, namely reversible maps defined by rational functions with 
coefficients in $K$.
If the field $K=\FF_q$ is a {\it finite field\/} with $q$ elements
---where $q$ is the power of a prime number $p$---
we are led to the study of reversibility on a phase space 
with $N=q^n$ points.
(For background reference on finite fields, see \cite{LidlNiederreiter}.)
More meaningful is the case in which $K$ is an algebraic number field,
e.g., $K=\QQ$. 
The same map $L$ may then be reduced over infinitely many distinct 
finite fields, of increasing size. 

We are concerned here with the case where the reductions of the component involutions $G$
and $H$ are involutory permutations of the phase space of $N=q^n$ points, so
that the reduction of their composition $L$ yields a permutation of this space. In reduction,
the fixed sets  $Fix\,G$ 
and $Fix\,H$ will now be finite.
Examples arising from dimension $n=2$ are the reversible planar polynomial automorphisms 
\cite{FriedlandMilnor,BaakeRoberts:05,GomezMeiss2}
which include the well-known area-preserving H\'enon map $L_1=H_1\circ G_1$,
\begin{equation}\label{eq:Henon}
G_1:(x,y)\mapsto (y,x)
\hskip 40pt
H_1:(x,y)\mapsto(x,-y+x^2+a)\qquad a\in K.
\end{equation}
If $K=\FF_p$, we find 
$$
\#Fix\,G_1\,=\, \#Fix\,H_1=p.
$$

Reversible polynomial automorphisms in any dimension provide other
examples. But certain rational maps in $n$ dimensions can induce 
permutations, e.g., $L_2=H_2\circ G_2$ of \cite{RobertsLamb}
\begin{eqnarray*}\label{eq:3DMap}
G_2:(x,y,z)&\mapsto&(x+e(2y-k)(z+e(y-k)),k-y,z+e(2y-k))\qquad e,k\in
K\\
H_2:(x,y,z)&\mapsto&(y\,(2-2x/F+x^2/F^2),x/F, -z) \qquad  F=1+(1-y)^2.
\end{eqnarray*}
If $K=\FF_p$, with $p \equiv 3 \mod 4$, the denominator in $H_2$ is
non-vanishing and
$$
\#Fix\,G_2\,=p^2,\qquad \#Fix\,H_2=p.
$$

The representability of $G$ and $H$ over a finite field $\FF_q$ imposes some 
restrictions on the parameters present in the involutions. 
For instance, if in \eqref{eq:Henon}, 
we take $a\in\QQ$, say, we are led to consider the set $P_a$ of the prime 
divisors of the denominator of $a$. Then the involution $H_1$ is representable 
over any finite field whose characteristic is not in $P_a$, such as the field
of integers modulo a prime $p\not\in P_a$. 
Additionally, one must verify that $G$ and $H$ remain involutions when reduced
to a finite field. This is an instance of the property of {\it good reduction} 
(for a discussion of reduction in arithmetic dynamics, see 
\cite[chapter 2]{Silverman}).

Representability assumed, the reduction of an algebraic map can be associated with 
interesting asymptotic phenomena, concerning the statistical behaviour 
of the periodic orbits  in the large field limit.
A main theme in this area of research is to ascertain whether such phenomena are 
{\it universal,} namely independent from the map, depending only on its structural 
properties e.g., reversibility.  Furthermore, in reduction,
the original dimension $n$ of the underlying algebraic map ceases to be important
as the finite phase space does not inherit the topology of the continuum $\RR^n$
or $\CC^n$.

If the reduction of a reversible algebraic map is a permutation,
it has been shown \cite{RobertsVivaldi:05} that over finite fields, 
reversibility manifests itself combinatorially. In particular,
the number of symmetric periodic orbits  ---those invariant under $G$--- 
equals 
\begin{equation} \label{eq:NumSym}
({\#Fix\,G+\# Fix\,H})/2.
\end{equation}
Strong experimental evidence \cite{RobertsVivaldi:05,Jogia}
suggests that the symmetric periodic orbits dominate the statistics over the 
asymmetric ones\footnote{This is in sharp contrast 
with the case of maps with real or complex coordinates \cite{MacKay}.}, 
so that \eqref{eq:NumSym} asymptotically counts the number of periodic orbits
in the permutation.
Furthermore, numerical experiments  \cite{RobertsVivaldi:05,Jogia} 
suggest there exists an asymptotic (large fields)
distribution of the periods of the orbits, given by
\begin{equation}\label{eq:R}
\cR(x)=1-e^{-x}(1+x)\qquad x\geq 0.
\end{equation}
(In these experiments, `large fields' refers to the limit $\FF_p,\,p\to\infty$;
there are other asymptotic regimes, see remark (3) in section 
\ref{section:ConcludingRemarks}.) 
The distribution (\ref{eq:R}) describes the limiting probability of the 
set of points belonging to cycles with scaled period not exceeding $x$.
This distribution is believed to be universal, within the class of 
reversible maps with a single time-reversal symmetry.
By contrast, maps which are not reversible appear to behave like random 
permutations \cite{RobertsVivaldi:05}; the associated distribution, 
also conjectured to be universal, is markedly different from that of reversible maps.
\footnote{The inadequacy of random permutations as a model for reversible 
maps on a discrete space was noted as far back as \cite{Rannou}.}
Finally, there are asymptotic period distributions for integrable systems 
\cite{RobertsVivaldi,JogiaRobertsVivaldi}. Even though these distributions
are map-dependent, they all feature a distinctive `quantization' of 
periods. With appropriate scaling, the allowed periods occur at
the reciprocals of the natural numbers, resulting in step-like 
distribution functions.
These various asymptotic phenomena have led to the development of
simple and effective tests for detecting integrability and reversibility 
in algebraic mappings \cite{RobertsVivaldi,RobertsVivaldi:05,JogiaRobertsVivaldi}.

In this paper we propose a combinatorial model for the orbits,
over a finite field, of a map with a single family of reversing
symmetries. Our model consists of the composition of 
two random involutions. 
The size $N$ of the space and the numbers $\#Fix\,G$ and $\#Fix\,H$ 
of fixed points of each involution are the parameters of the model.
By letting $N$ increase, we build a sequence of probability spaces. 
Under very general conditions, we derive the distribution (\ref{eq:R})
found in \cite{RobertsVivaldi:05,Jogia}.

We prove the following theorem, which appears in section
\ref{section:RandomInvolutions} as theorems 
\ref{theorem:Distribution} and \ref{theorem:AsymmetricCycles}.

\bigskip
\noindent{\bf Theorem A.} {\sl 
Let $(G,H)$ be a pair of random involutions of a set $\Omega$ with $N$ 
points, and let $g=\#Fix\,G$ and $h=\#Fix\,H$. Let $\cR_N(x)$
be the expectation value of the portion of $\Omega$ occupied by cycles 
of $H\circ G$ with period less than $2xN/(g+h)$, computed with
respect to the uniform probability. 
If, with increasing $N$, $g$ and $h$ satisfy the conditions
\begin{equation}\label{eq:ghAsymptotics}
\lim_{N\to\infty} g(N)+h(N)=\infty\qquad
\lim_{N\to\infty}\frac{g(N)+h(N)}{N}=0
\end{equation}
then for all $x\geq 0$, we have the limit $\cR_N(x)\to\cR(x)$, 
the distribution function in (\ref{eq:R}). 
Moreover, almost all points in $\Omega$ belong to symmetric cycles.
}

\bigskip\noindent
This result is formulated in terms of the sequence $g+h$, and not of the 
two sequences individually.  The conditions on the growth rate of $g+h$
are quite mild, hence the wide scope of applicability of the theorem.
If $g$ and $h$ do not grow at the same rate, an interesting dynamical
phenomenon occurs: asymptotically, almost all periodic orbits have even 
period (corollary \ref{corollary:EvenOnly}).
For $n$-dimensional maps over the finite field $\FF_q$, one has $N=q^n$,
while $g$ and $h$ will typically grow algebraically (see remark 1 of section 
\ref{section:ConcludingRemarks}).
As mentioned above, increasing the dimension does not introduce new difficulties, 
because the absence of a significant topological structure 
makes the combinatorial model effectively dimension-independent. 
All that is needed is to keep track of the possible cardinalities of 
the fixed sets of the two involutions.

In \cite{RobertsVivaldi:05}, 
experimental evidence from reductions of the H\'enon map \eqref{eq:Henon}
suggested that the probability of occurrence
of several cycles with the same period followed a Poisson law,
with parameter depending exponentially on the (rescaled) period.
These findings motivate the study of the statistics of the occurrence 
of repeated periods in the random involutions model.
In section \ref{section:Repetitions} we consider this question,
and prove the following result.

\vspace*{10pt}
\noindent{\bf Theorem B.} {\sl 
With the notation of theorem A, let $\mu(t,i)$ be the probability that
$(G,H)$ has $i$ cycles of period $t$. Assume that the sequence
\begin{equation}\label{eq:f}
f(N)=\begin{cases} \displaystyle \frac{gh}{N} & \mbox{if $t$ is odd}\\
\noalign{\vspace*{3pt}}
\displaystyle \frac{g^2+h^2}{2N} & \mbox{if $t$ is even}
\end{cases}
\end{equation}
has the (possibly infinite) limit 
$$
c=\lim_{N\to\infty} f(N)
$$
and let
\begin{equation}\label{eq:y}
y=\frac{(t-1)(g+h)}{2N}-\ln(f).
\end{equation}
Then, if $c\not=0$ we have, as $N\to\infty$
\begin{equation}\label{eq:muAsymptotic}
\mu(t,i)\sim 
e^{\displaystyle -\alpha}\displaystyle \frac{\alpha^i}{i!}\qquad 
\alpha=\begin{cases}
c & \mbox{if $t$ is constant}\\
e^{-y(t)} & \mbox{if $y$ is constant}
\end{cases}
\end{equation}
while if $c=0$ we have $\mu(t,i)\sim\delta_i$, where $\delta$ is Kronecker's delta.
}
\vspace*{10pt}

As it will appear in the proof of theorem B, in the case of constant $y$, 
for convergence to a Poisson distribution it is only required that 
$f(N)\not\to 0$, not the actual existence of a non-zero limit $c$. We conclude
section 3 by showing the agreement between the prediction of theorem B
with the actual repetitions $\mu(t,i)$ for a small $t$ in the H\'enon map.
The latter 
are governed by the roots $\mod p$ of certain polynomial functions of one variable,
with a connection to the Galois groups of these polynomials.

In section \ref{section:ConcludingRemarks} we discuss the question
of the asymptotic (large field) growth rate of the parameters
associated with a rational map over a finite field, namely the 
numbers $g$ and $h$ of fixed points of the involutions.
In this case $N$ is the power of a prime number.
We show that the growth conditions of theorem A on $g+h$ are 
easily satisfied, while $f(N)$ in theorem B depends algebraically
on $N$, and all scenarios considered in the theorem can be realized.
Finally, we also sketch how the random involution model might be
extended to the case of reversible rational maps with singularities.

\section{Composition of random involutions}\label{section:RandomInvolutions}

Let $N$ be a positive integer, and let $g,h$ be integers in the
range $1\leq g,h\leq N$ such that $N-g$ and $N-h$ are both even.
We consider ordered pairs $(G, H)$ of random involutions on a set $\Omega$
of $N$ points, which fix $g$ and $h$ points, respectively. Thus $g=\#Fix\,G$, 
$h=\#Fix\,H$, and the cycle decomposition of $G$ consists of $g$ 
fixed points and $(N-g)/2$ two-cycles, and similarly for $H$. 
For each pair $(G,H)$ of involutions we consider the composition 
$L=H\circ G$, which is a reversible permutation of $\Omega$. 

Let ${\bf G}$ and ${\bf H}$ be the sets of all involutions with 
the given parameters. We regard the space 
\begin{equation}\label{eq:EGH}
{\bf E}={\bf E}(g,h,N)={\bf G}\times {\bf H}
\end{equation} 
as a probability space, with the uniform probability. 
We define the random variable $P_t:{\bf E}\to \RR$ to be the fraction 
of the space $\Omega$ occupied by the $t$-cycles of the map $H\circ G$, namely
\begin{equation}\label{eq:RandomVariable}
P_t=\frac{1}{N}\,\#\{x\in \Omega\,:\, x 
    \mbox{ has minimal period $t$ under } H\circ G\}\qquad t=1,2,\ldots .
\end{equation}
For given $N$, the sequence $(P_t)$ has only a finite number of non-zero terms.

The main object under study is the distribution function 
\begin{equation}\label{eq:Rn}
\cR_N(x)=\sum_{t=1}^{\lfloor xz\rfloor }\langle P_t\rangle
\hskip 30pt x\geq 0
\end{equation}
where $\lfloor\cdot\rfloor$ is the floor (greatest integer part) function, 
$z$ is a scaling parameter to be specified below, and the average 
$\langle\cdot\rangle$ is computed with respect to the uniform 
probability on ${\bf E}$ (i.e., $\langle P_t\rangle=(\# {\bf E})^{-1}\sum_{\bf E} P_t$). 
If $\lfloor xz\rfloor=0$, the sum is empty, and $\cR_N(x)$ is defined to be zero. 

The asymptotic parameter is $N$. We shall assume that $g=g(N)$ 
and $h=h(N)$ are positive integer sequences satisfying \eqref{eq:ghAsymptotics}.
These sequences  determine the sequence of probability 
spaces to be considered.
To get non trivial asymptotics, the scaling parameter $z$ must grow 
with $N$ at an appropriate rate.
To this end, we define the rational sequence
\begin{equation}\label{eq:z}
z(N)=\frac{2N}{g(N)+h(N)}
\end{equation}
which, due to (\ref{eq:ghAsymptotics}), diverges to infinity.
The quantity $z(N)$ is the ratio of the size of the phase space
to the number of symmetric cycles ---cf.~equation (\ref{eq:NumSym}).
Its full significance will be clarified  at the end of 
this section (corollary \ref{corollary:z}).

The main purpose of this section is to prove theorem A of
the introduction, which is an amalgamation of theorems
\ref{theorem:Distribution} and \ref{theorem:AsymmetricCycles}.
\begin{theorem}\label{theorem:Distribution}
The distribution function $\cR_N$ defined in (\ref{eq:Rn}), with scaling 
sequence $z(N)$ given by (\ref{eq:z}), admits the following limit
$$
\lim_{N\to\infty}\cR_N(x)=\cR(x)=1-e^{-x}(1+x)\hskip 30pt x\geq 0.
$$
\end{theorem}

\medskip
It turns out that to establish theorem \ref{theorem:Distribution}
it suffices to consider symmetric orbits only, because they provide 
the dominant contribution to the probability $\langle P_t\rangle$ 
in equation (\ref{eq:Rn}) ---see theorem \ref{theorem:AsymmetricCycles} below.
Accordingly, we separate out the contribution deriving from symmetric
and asymmetric cycles as follows
\begin{equation}\label{eq:Decomposition}
P_t=P_t^{(s)}+P_t^{(a)}.
\end{equation}

\begin{figure}[h]
\hspace*{80pt}\epsfig{file=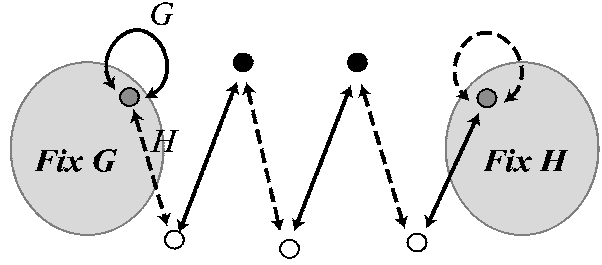,width=20cm,height=4cm}\hfil
\caption{\label{fig:Arc}\rm\small Construction of a 4-arc
corresponding to a 7-cycle. 
The solid and dashed lines denote the action of $G$ and $H$, 
respectively. The orbit of $L=H\circ G$ from $Fix\, G$ to $Fix\, H$
is represented by white circles; the black circles represent 
the return orbit from $Fix\, H$ to $Fix\, G$.}
\end{figure}

We begin with symmetric cycles of odd period $t=2k-1$. It is known 
\cite{Devaney,DeVog,LambRoberts,MacKay,RobertsQuispel}
that a cycle is of this type if and only if it corresponds to an orbit of 
$L$ starting from $Fix\,G$ and reaching $Fix\,H$ for the first time after 
$k$ iterations.
Here and below, we follow the cycle dynamics keeping track of the action of 
$G$ and $H$ separately, which leads to a maximal sequence of $2k$ points
\begin{equation}\label{eq:Arc}
(x_1,y_1,\ldots,x_{k},y_k)
\qquad y_j=G(x_j),\quad x_{j+1}=H(y_j)
\quad k\geq 1, \quad j=1,\ldots,k
\end{equation}
with the property that all $x_j$ are distinct.
We call such a sequence a {\it $k$-arc}.
In the case of odd period, we have the condition
\begin{equation}\label{eq:ArcOdd}
x_1\in Fix\,G,\qquad y_{k}\in Fix\,H
\end{equation}
(hence $y_1=x_1$) together with the requirement that the sequence 
contains $t=2k-1$ distinct points (see figure \ref{fig:Arc}).
We must now determine the total number of $k$-arcs of this type
generated by the elements of ${\bf E}$.

Throughout the paper, the underlined exponent will denote the 
falling factorial powers
\begin{equation}\label{eq:FallingFactorial}
n^{\underline{0}}=1
\hskip 40pt
n^{\underline{a}}=n(n-1)(n-2)\cdots (n-a+1),\qquad a\geq 1.
\end{equation}
The information we require is contained in the following two lemmas.

\begin{lemma}\label{lemma:E} 
Let ${\bf E}={\bf E}(g,h,N)$ be as in equation (\ref{eq:EGH}). Then
\begin{equation}\label{eq:E}
\#{\bf E}(g,h,N)=\frac{1}{2^{N-(g+h)/2}}\,
\frac{(N!)^2}{g!\,h!\,\left((N-g)/2\right)!\left((N-h)/2\right)!}.
\end{equation}
\end{lemma}

\bigskip

\proof
The cardinality of the set of all involutions on $N$ points, 
with $g$ fixed points, is given by
$$
{N \choose g}\,\frac{(N-g)!}{2^{(N-g)/2}((N-g)/2)!}
$$
where the first term gives the number of ways of choosing the
fixed points, and the second counts the number of ways of 
partitioning the remaining points into unordered pairs.
Multiplying the expression above with the corresponding expression 
for $h$ fixed points, we obtain
$$
\#{\bf E}(g,h,N)=\frac{1}{2^{N-(g+h)/2}}\frac{(N!)^2}{g!\, h!\, ((N-g)/2)!\,((N-h)/2)!}
$$
as desired. 
\endproof

\begin{lemma}\label{lemma:POdd} Let $P_t^{(s)}$ be as above.
The following holds
\begin{equation}\label{eq:POdd}
\langle P_{2k-1}^{(s)}\rangle=\frac{2k-1}{N}\,N^{\underline{2k-1}}\,\,\,\cE_{2k-1}\, 
\qquad k\geq 1
\end{equation}
where 
\begin{eqnarray}
\cE_{2k-1}
&=&\frac{\#{\bf E}(g-1,h-1,N-2k+1)}{\#{\bf E}(g,h,N)}\nonumber\\
&=&\begin{cases}
\displaystyle\frac{gh}{N^2} & \mbox{if } k=1\\
\displaystyle \frac{gh}{(N^{\underline{2k-1}})^2}\,
\prod_{j=0}^{k-2}(N-g-2j)(N-h-2j)
        & \mbox{if } k>1.\end{cases} \label{eq:cEOdd}
\end{eqnarray}
\end{lemma}

\bigskip

\proof Let $(x_1,y_1,\ldots,x_{k},y_k)$ be the required $k$-arc
(see equation (\ref{eq:Arc}) and figure \ref{fig:Arc}).
The initial point $x_1$ belongs to $Fix\,G$, so that $y_1=x_1$.
The first half of the cycle consists of the points $x_1,\ldots,x_k$, 
the second half of the points $y_k,\ldots,y_2$, 
giving $2k-1$ distinct points in total.
Hence the desired number of arcs is
$$
N^{\underline{2k-1}}.
$$
Each $k$-arc occurs with multiplicity determined by the unconstrained 
action of the pairs of involutions on the rest of the space, which is
given by
$$
\#{\bf E}(g-1,h-1,N-2k+1).
$$
Thus the average value of $P_{2k-1}^{(s)}$ is given by
$$
\langle P_{2k-1}^{(s)} \rangle=\frac{2k-1}{N}\,N^{\underline{2k-1}}\,
\frac{\#{\bf E}(g-1,h-1,N-2k+1)}{\#{\bf E}(g,h,N)}.
$$
Finally, the product formula for $\cE_{2k-1}$ is obtained from lemma 
\ref{lemma:E} with a straightforward computation.
\endproof


By analogy with equation (\ref{eq:Rn}), we define the distribution function 
$\cR^{(s,o)}_N$ for odd symmetric cycles as
\begin{equation}\label{eq:RnOdd}
\cR_N^{(s,o)}(x)=\sum_{k=1}^{\lfloor (xz+1)/2\rfloor}\langle P_{2k-1}^{(s)}\rangle
\qquad x\geq 0.
\end{equation}
We have
\begin{theorem}\label{theorem:DistributionOdd}
For any sequences $g,h$ satisfying (\ref{eq:ghAsymptotics}), 
with $z$ as in (\ref{eq:z}), we have,
as $N\to\infty$
$$
\cR_N^{(s,o)}(x)\sim\frac{2gh}{(g+h)^2}\,\cR(x)
\hskip 30pt x\geq 0.
$$
\end{theorem}
Note that the limit distribution is symmetrical in $g$ and $h$.
To prove theorem \ref{theorem:DistributionOdd} we need a lemma.
\begin{lemma}\label{lemma:Basic}
Let $g(N),h(N)$ be as in (\ref{eq:ghAsymptotics}), let
$m=m(N)$ be a positive integer sequence, and let
$$
\lambda=\lambda(N)=\left(1-\frac{g}{N}\right)\left(1-\frac{h}{N}\right)
\hskip 40pt
S_m=\sum_{k=1}^{m}(2k-1)\lambda^{k-1}.
$$
Then, as $N\to\infty$
\begin{equation}\label{eq:Basic}
S_m\sim \frac{2}{\kappa^2}\bigl[1-e^{-m\kappa}(m\kappa+1)\bigr]
\hskip 30pt\kappa=\frac{g+h}{N}.
\end{equation}
\end{lemma}

\proof
For fixed $N$, $m$ and $\lambda$ are fixed, and the quantity $S_m$ is the sum 
of an arithmetico-geometric progression, with value
$$
S_m=\frac{1}{(1-\lambda)^2}
\Bigl\{1+\lambda-\lambda^m \bigl[2m(1-\lambda)+1+\lambda\bigr]\Bigr\}.
$$
As $N\to\infty$, we have, from (\ref{eq:ghAsymptotics})
$$
\lambda\sim 1\qquad  
(1-\lambda)\sim \kappa\qquad
\ln(\lambda)\sim -\kappa 
$$
so we find
$$
S_m\sim \frac{1}{\kappa^2}\bigl[2-e^{m\ln(\lambda)}(2m\kappa+2)\bigr]
\sim \frac{2}{\kappa^2}\bigl[1-e^{-m\kappa}(m\kappa+1)\bigr]
$$
as desired.
\endproof

\bigskip 
\noindent {\it Proof of theorem \ref{theorem:DistributionOdd}.} \/
From lemma \ref{lemma:POdd}, we find, for $k\geq 1$
\begin{eqnarray}
\langle P_{2k-1}^{(s)}\rangle&=&\frac{gh(2k-1)}{N}\,
\frac{\prod_{j=0}^{k-2}(N-g-2j)(N-h-2j)}
{N^{\underline{2k-1}}}\nonumber\\
&=&\frac{gh(2k-1)}{N(N-(2k-2))}\,
\prod_{j=0}^{k-2}\left(1-\frac{g-1}{N-(2j+1)}\right)\left(1-\frac{h}{N-2j}\right).
 \label{eq:AvgpOdd}
\end{eqnarray}

Now fix $m$, with $0<2m<N$, and define
\begin{equation}\label{eq:LambdaPlusMinus}
\lambda_-= \left(1-\frac{g}{N-m}\right)\left(1-\frac{h}{N-m}\right)
\qquad
\lambda_+=\left(1-\frac{g-1}{N}\right)\left(1-\frac{h}{N}\right).
\end{equation}
For all $k$ such that $1\leq 2k-2\leq 2m$, we have the inequalities
\begin{equation}\label{eq:Bounds}
\lambda_-^{k-1}
<\prod_{j=0}^{k-2}\left(1-\frac{g-1}{N-(2j+1)}\right)\left(1-\frac{h}{N-2j}\right)
<\lambda_+^{k-1}
\end{equation}
which, together with (\ref{eq:AvgpOdd}), give
\begin{equation} \label{eq:Compare}
\frac{gh (2k-1)}{N^2}\,\,\lambda_-^{k-1}
<\langle P_{2k-1}^{(s)}\rangle
<\frac{gh(2k-1)}{N(N-m)}\,\,\lambda_+^{k-1}.
\end{equation}
Summing over $k$, we obtain the bounds
\begin{equation}\label{eq:RBounds}
\cR_-\leq \sum_{k=1}^m \langle P_{2k-1}^{(s)}\rangle \leq\cR_+, 
\end{equation}
where
$$
\cR_-=\frac{gh}{N^2}\sum_{k=1}^{m}(2k-1)\lambda_-^{k-1}
\hskip 40pt
\cR_+=\frac{gh}{N(N-m)}\sum_{k=1}^{m}(2k-1)\lambda_+^{k-1}.
$$
Let now $m=m(N)$ be a positive integer sequence such that
$$
\lim_{N\to\infty}m(N)=\infty
\hskip 40pt
m=o(N).
$$
Then $\lambda_\pm\sim\lambda$ and 
applying lemma \ref{lemma:Basic}, we obtain, as $N\to\infty$
\begin{equation}\label{eq:Asymptotics}
\cR_-\sim\cR_+\sim 
\frac{gh}{N^2}S_{m}\sim
\frac{2gh}{(g+h)^2}\bigl[1-e^{-m\kappa}(m\kappa+1)\bigr].
\end{equation}

We now fix a non-negative real number $x$. 
For the convergence of the expression above, 
we specialize the sequence $m(N)$ as follows
\begin{equation}\label{eq:m}
m\kappa=x\qquad \mbox{or}\qquad m(N)=\lfloor xz(N)/2\rfloor
\end{equation}
where $\kappa$ and $z$ were defined in equations (\ref{eq:Basic}) and 
(\ref{eq:z}), respectively.
Due to (\ref{eq:ghAsymptotics}), $m(N)=o(N)$, and for all 
sufficiently large $N$, we have $0<m(N)<N/2$.
Hence (\ref{eq:RBounds}) is valid, and our result now follows from 
(\ref{eq:RnOdd})  and (\ref{eq:Asymptotics}), noting that 
$(xz(N)+1)/2\sim xz(N)/2$.
\endproof

The sequence $m(N)$ given in (\ref{eq:m}) could be defined by the
more general requirement that the sequence $m(N)\kappa(N)$ converge
to a positive real number $c$, where the limiting cases $c=0,\infty$ 
are excluded as they would lead to trivial distribution functions.
We have set $c=1$, for normalization.
We remark that the choice of this constant does not affect the 
coefficient of the distribution function, which determines its 
limiting value for large arguments. This fact will become relevant 
in the proof of theorem \ref{theorem:Distribution} below.

\medskip
We now turn to even symmetric cycles, of period $t=2k$. Each cycle of this 
type corresponds to an orbit originating on $Fix\,G$ (or $Fix\,H$), and
returning to it for the first time after $k$ iterations of $L$.
The arc associated to such a cycle (see equation (\ref{eq:Arc})) is
a $(k+1)$-arc, with the property that $x_1$ and $x_{k+1}$ are in $Fix\,G$.
Therefore $x_1=y_1$ and $x_{k+1}=y_{k+1}$, and the total number of
arcs of this type is equal to $N^{\underline{2k}}/2$, where the denominator 
2 prevents arcs beginning and ending on ${\rm Fix}\, G$ from being counted twice.
The denominator 2 cancels with the 2 in the period $2k$, so we obtain
\begin{equation}\label{eq:PEven}
\langle P_{2k,g}^{(s)}\rangle=\frac{k}{N}\,N^{\underline{2k}} \,\cE_{2k}
\qquad k\geq 1
\end{equation}
where
\begin{eqnarray}
\cE_{2k}
&=&\frac{\#{\bf E}(g-2,h,N-2k)}{\#{\bf E}(g,h,N)}\nonumber\\
&=&\begin{cases} 
    \displaystyle \frac{g(g-1)(N-h)}{N(N-1)^2}& \mbox{if } k=1\\
    \displaystyle  \frac{g(g-1)}{N^{\underline{2k}}} \,
    \prod_{j=0}^{k-2}(N-g-2j)\prod_{j=0}^{k-1}(N-h-2j) & \mbox{if } k>1.
    \end{cases}\label{eq:cEEven}
\end{eqnarray}

Putting together equation (\ref{eq:PEven}) and (\ref{eq:cEEven}), we obtain,
for $k\geq 2$
$$
\langle P_{2k,g}^{(s)}\rangle=\frac{g(g-1)k}{N(N-(2k-1))}\,
\prod_{j=0}^{k-2}\left(1-\frac{g-1}{N-(2j+1)}\right)
\prod_{j=0}^{k-1}\left(1-\frac{h}{N-2j}\right)
$$
to be compared with equation (\ref{eq:AvgpOdd}).

Exchanging $g$ by $h$ in the formulae above gives the twin quantity 
$P_{2k,h}^{(s)}$ for $2k$-cycles with two points of $Fix\,H$. 
The distribution function $\cR^{(s,e)}_N$ for even symmetric cycles 
is obtained by adding the contributions from the two fixed sets
\begin{equation}\label{eq:RnEven}
\cR_N^{(s,e)}(x)=
\sum_{k=0}^{\lfloor xz/2\rfloor}\,
 \bigl(\langle P_{2k,g}^{(s)}\rangle+\langle P_{2k,h}^{(s)}\rangle\bigr)
\qquad x\geq 0.
\end{equation}

An analysis very similar to that of theorem \ref{theorem:DistributionOdd}
leads to the following result, whose proof we omit for the sake of brevity.
\begin{theorem}\label{theorem:DistributionEven}
For any sequences $g,h$ satisfying (\ref{eq:ghAsymptotics}), 
with $z$ as in (\ref{eq:z}), we have,
as $N\to\infty$
$$
\cR_N^{(s,e)}(x)\sim\frac{g^2+h^2}{(g+h)^2}\,\cR(x)
\hskip 30pt x\geq 0.
$$
\end{theorem}

\medskip

\noindent {\sc Proof of theorem \ref{theorem:Distribution}.} \/
From theorems \ref{theorem:DistributionOdd} and \ref{theorem:DistributionEven}
we deduce that the probability distribution associated to symmetric cycles
is given by
$$
\cR_N^{(s,o)}(x)+\cR_N^{(s,e)}(x) \sim \frac{2gh+g^2+h^2}{(g+h)^2}\,\cR(x)=\cR(x).
$$
Because
$$
\lim_{x\to\infty}\cR(x)=1
$$
this distribution function is properly normalized (cf.~ remarks
following the proof of theorem \ref{theorem:DistributionOdd}), 
and so it accounts for all cycles.
\endproof

It is implicit from theorem \ref{theorem:Distribution} that, asymptotically, 
asymmetric cycles have zero probability. Asymmetric $t$-cycles of 
$L=H \circ G$ come in pairs, mapped to one another under $G$.
They correspond to $t$-arcs of the type
\begin{equation} \label{eq:ArcAsym}
(x_1,y_1,\ldots,x_{t},y_t) \quad
y_j=G(x_j)\ne x_j,\; \; x_{j+1}=H(y_j)\ne y_j,
\quad j=1,\ldots,t
\end{equation}
i.e., they have $2t$ distinct entries. 
We can now determine the decay rate of the probability
of asymmetric cycles.

\begin{theorem} \label{theorem:AsymmetricCycles}
Consider a composition of random involutions on a set of $N$
points, whose fixed sets satisfy (\ref{eq:ghAsymptotics}). 
As $N\to\infty$, the measure of the set of points that belong 
to asymmetric cycles is asymptotic to $(g(N)+h(N))^{-1}$;
in particular, almost all cycles are symmetric.
\end{theorem}

\bigskip
\proof 
Recall from (\ref{eq:Decomposition}) that $P_t^{(a)}$ is the 
probability of finding an asymmetric $t$-cycle. We compute 
its average over ${\bf E}$:
\begin{eqnarray}
\langle P_t^{(a)}\rangle
&=&\frac{N^{\underline{2t}}}{N}\,
\frac{\#{\bf E}(g,h,N-2t)}{\#{\bf E}(g,h,N)}\nonumber\\
&=&\frac{1}{N}\,\prod_{j=0}^{t-1}\left(1-\frac{g-1}{N-(2j+1)}\right)\left(1-\frac{h}{N-2j}\right).
\end{eqnarray}
From the above discussion, the right-hand side of this equation involves
counting $t$-arcs like \eqref{eq:ArcAsym}. Note that the period of the
orbit is not present as a factor. This is because it is absorbed by the 
possibility to create $t$-arcs that are one and the same from 
\eqref{eq:ArcAsym} by cyclic permutations of the  
$\{x_j,y_j\}$ pairs ($t$ possibilities) and the switch 
$x_j \leftrightarrow y_j$ together with reversing the order 
of the switched pairs ($2$ possibilities).

Let $\cR_N^{(a)}$ be the associated distribution function
$$
\cR_N^{(a)}(x)=\sum_{t=1}^{\lfloor xz\rfloor} \langle P_t^{(a)}\rangle.
$$
We now proceed as in the proof of theorem, with $\lambda_\pm$ and
$z$ given by (\ref{eq:LambdaPlusMinus}) and (\ref{eq:z}), 
respectively. Summing the expression above over $t$ in the
appropriate range gives the bounds
$$
\cR_-\leq \cR_N^{(a)}(x) \leq\cR_+, 
$$
where
$$
\cR_-=\frac{1}{N}\sum_{t=1}^{\lfloor xz(N)\rfloor}\lambda_-^{t}
\hskip 40pt
\cR_+=\frac{1}{N}\sum_{t=1}^{\lfloor xz(N)\rfloor}\lambda_+^{t}.
$$
This time the bounds are geometric progressions;
as $N\to \infty$, we find
$$
\cR_-\sim\cR_+\sim\frac{1-e^{-x}}{g(N)+h(N)},
$$
as desired.
\endproof

The theorems of this section have some immediate consequences.
Consider the sequence $z(N)$, defined in (\ref{eq:z}). 
Its denominator $g(N)+h(N)$ is twice the number of symmetric cycles
\cite[proposition 1]{RobertsVivaldi:05}, and from theorem 
\ref{theorem:Distribution}, asymptotically, the symmetric cycles have full measure.
These considerations give the following characterization of $z(N)$.

\begin{corollary} \label{corollary:z}
Asymptotically, the sequence $z(N)$ defined in (\ref{eq:z}) 
represents the average cycle length of the map $H\circ G$.
\end{corollary}

An immediate consequence of theorems \ref{theorem:DistributionOdd} 
and \ref{theorem:DistributionEven} is the following.

\begin{corollary} \label{corollary:EvenOnly}
Let $g,h$ be as above. If $g(N)=o(h(N))$ or $h(N)=o(g(N))$, 
then almost all cycles of $H\circ G$ have even period.
\end{corollary}

\section{Occurrence of repeated periods}\label{section:Repetitions}

This section first gives the proof of theorem B, stated in the introduction.
The asymptotic relations (\ref{eq:ghAsymptotics}) will be assumed to hold, 
together with the condition $h=O(g)$. The latter does not entail loss of 
generality, as we may always exchange $g$ and $h$, which amounts to considering
the inverse mapping $L^{-1}$.  After the proof, we compare the results
of theorem B with the mechanism that determines the repetition statistics for
small period in the H\'enon map \eqref{eq:Henon}.

\medskip

\noindent {\sc Proof of theorem B.} \/
Let $\mu(t,i)$ be the probability that a randomly chosen 
$(G,H)\in {\bf E}$ has exactly $i$ cycles of period $t$.
We wish to determine the asymptotic behaviour of $\mu$;
due to theorem \ref{theorem:AsymmetricCycles}, it will 
suffice to consider symmetric cycles only.

We introduce some notation. Let $\gamma$ be a $t$-cycle in $\Omega$,
and let $\Gamma$ be any set of disjoint $t$-cycles.
We stipulate that the first point of each cycle belongs to some 
symmetry line; then there are
\begin{equation}\label{eq:ChoosingCycles}
\frac{N^{\underline{it}}}{i!}
\end{equation}
ways of choosing such cycles.
We define
\begin{equation}\label{eq:SetsE}
E_\gamma=\{(G,H)\in {\bf E}\,:\, \gamma \mbox{ is a $t$-cycle of } (G,H)\}
\hskip 40pt
E_\Gamma=\bigcap_{\gamma\in\Gamma}\, E_\gamma.
\end{equation}
Note that if the cycles of $\Gamma$ were not disjoint, then $E_\Gamma$ 
would be empty.

Given $\Gamma$, with $\#\Gamma=i$, we are interested in the number
of pairs in {\bf E} that support the cycles of $\Gamma$, and that do 
not have any $t$-cycle in the complement $\Omega\setminus\Gamma$.
The parameters of the involutions acting on $\Omega\setminus\Gamma$ 
are determined by the type of $t$-cycles being considered (odd period 
symmetric, even period symmetric on $Fix\,G$, etc.).
Accordingly, we'll consider the following sets
of pairs of involutions

\medskip

\begin{equation}\label{eq:Es}
\hfil\vcenter{\halign{
$#$\hfil\qquad&#\hfil\cr
{\bf E}(g-i,h-i,N-it)&symmetric, odd period\cr
\noalign{\vskip 1pt}
{\bf E}(g-2i,h,N-it)&symmetric, even period, on $Fix\,G$\cr
\noalign{\vskip 1pt}
{\bf E}(g,h-2i,N-it)&symmetric, even period, on $Fix\,H$\cr
}}\hfil
\end{equation}

\medskip

\noindent
whose cardinality is given by lemma \ref{lemma:E}.
Whenever we don't require specialization, we'll denote any of the above sets 
by the common symbol ${\bf E}_{i}$.
The pairs of involutions in ${\bf E}_{i}$ act on the space $\Omega_i$ 
with $N-it$ elements.

Let ${\bf \bar E}_i$ be the set of elements of ${\bf E}_{i}$ that 
have no $t$-cycles at all. By the inclusion-exclusion principle, we find
\begin{eqnarray}
\#{\bf \bar E}_i&=&\#\left({\bf E}_{i}\setminus \bigcup_{\gamma\in\Omega} 
    E_\gamma\right)= \sum_\Gamma (-1)^{\#\Gamma}\# E_\Gamma\nonumber\\
 &=&\sum_{n\geq 0}(-1)^n\sum_{\#\Gamma=n}\,\#E_\Gamma.\label{eq:nt}
\end{eqnarray}
The number of non-zero summands is clearly finite.

To compute the inner sum for fixed $t$ and $i$, we first choose $n$ distinct
cycles, and this can be done in $(N-it)^{\underline{nt}}/n!$ ways. 
Then we let the involutions act freely on the remaining space, 
which gives $\#{\bf E}_{i+n}$ possibilities.
Thus
$$
\sum_{\#\Gamma=n}\,\#E_\Gamma
=
\frac{(N-it)^{\underline{nt}}}{n!}\,\#{\bf E}_{i+n}.
$$

Considering the equations (\ref{eq:ChoosingCycles}) and (\ref{eq:nt}), we obtain
\begin{equation}\label{eq:mu}
\mu(t,i)= \frac{N^{\underline{it}}}{i!} \frac{\#{\bf \bar E}_i}{\#{\bf E}_0}
= \frac{N^{\underline{it}}}{i!}\,\sum_{n=0}^M(-1)^n
\frac{(N-it)^{\underline{nt}}}{n!}\frac{\#{\bf E}_{i+n}}{\#{\bf E}_0}
\end{equation}
where $M=\lfloor N/t-i\rfloor$.

We specialize formula (\ref{eq:mu}) to the case of odd period $t$, 
using the notation $\mu^{(o)}$. For this purpose, we consider the 
first expression ${\bf E}$ in (\ref{eq:Es}), to obtain
$$
\mu^{(o)}(i,t)=
\frac{N^{\underline{it}}}{i!}\,\sum_{n=0}^M (-1)^n
\frac{(N-it)^{\underline{nt}}}{n!}\,
\frac{\#{\bf E}(g-(i+n),h-(i+n),N-(i+n)t)}{\#{\bf E}(g,h,N)}.
$$
Inserting the values of $\#{\bf E}$ from lemma \ref{lemma:E}, we obtain,
after some manipulation
\begin{equation}\label{eq:muOdd}
\mu^{(o)}(t,i)=\frac{1}{i!}\sum_{n=0}^M\frac{(-1)^n}{n!}
  {\mathcal A}^{(o)}(i,n)\,{\mathcal B}^{(o)}(i,n)
\end{equation}
where
\begin{eqnarray}
{\mathcal A}^{(o)}(i,n)&=&
 \displaystyle \frac{g^{\underline{i+n}}\,h^{\underline{i+n}}}
 {(N-(i+n)(t-1))^{\underline{i+n}}} 
 \label{eq:Ao}\\
{\mathcal B}^{(o)}(i,n)&=& \begin{cases}
1 & \mbox{if $t=1$}\\
\displaystyle \prod_{j=0}^{(i+n)(t-1)/2-1}\left(1-\frac{g-1}{N-(2j+1)}\right)
          \left(1-\frac{h}{N-2j}\right) & \mbox{if $t\geq1$}.\end{cases}
\label{eq:Bo}
\end{eqnarray}

With a similar argument we derive the probability $\mu^{(e)}$ for even period
\begin{equation}\label{eq:muEven}
\mu^{(e)}(t,i)=\frac{1}{i!}\sum_{n=0}^M\frac{(-1)^n}{n!}
  {\mathcal A}^{(e)}(i,n)\,{\mathcal B}^{(e)}(i,n)
\end{equation}
where
\begin{eqnarray}
{\mathcal A}^{(e)}(i,n)&=&
\displaystyle\frac{1}{2}\,\frac{g^{\underline{2(i+n})}+h^{\underline{2(i+n})}}
{(N-(i+n)(t-1))^{\underline{i+n}}} 
\label{eq:Ae}\\
{\mathcal B}^{(e)}(i,n)&=&
\prod_{j=0}^{(i+n)(t-2)/2-1}\left(1-\frac{g-1}{N-(2j+1)}\right)
                           \left(1-\frac{h}{N-2j}\right)\label{eq:Be}\\
&&
\qquad \times 
\prod_{j=(i+n)(t-2)/2}^{(i+n)t/2-1}\left(1-\frac{h}{N-j}\right).\nonumber
\end{eqnarray}
As in the analysis of section \ref{section:RandomInvolutions}, this 
expression takes into account cycles on both $Fix\,G$ and $Fix\,H$, 
while the factor of 2 compensates for the double-counting that
results from the presence of two points on each symmetry line.
\bigskip

We now consider the asymptotic behaviour of $\mu(t,i)$, for fixed $i$. 
Let $y$ be given by (\ref{eq:y}), which we rewrite as
\begin{equation}\label{eq:yNew}
y=x-\ln(f)
\hskip 40pt
x(t)=\frac{(t-1)(g+h)}{2N}.
\end{equation}
For the purpose of scaling, the period $t$ will be allowed to vary with $N$.

Let us consider the expression
\begin{equation}\label{eq:Lambda}
\Pi(a,g,h)=\left(1-\frac{g}{N-a}\right)^{J/2}\left(1-\frac{h}{N-a}\right)^{J/2}
\qquad J=(i+n)(t-1)\qquad 0\leq a < N,
\end{equation}
which, for $a=o(N)$, provides the dominant contribution to
the products ${\mathcal B}(i,n)$, in both cases. 
A straightforward calculation gives
\begin{eqnarray}
\Pi(a,g,h)
 &=&e^{-(i+n)x}\left[1+O\left(\frac{nxa}{N}\right)
   +O\left(\frac{nxg}{N}\right)\right]\label{eq:BigO}.
\end{eqnarray}
In view of the bounds
$$
\Pi(J,g,h)\leq {\mathcal B}(i,n)\leq \Pi(0,g-1,h)
$$
we compute
\begin{eqnarray*}
\Pi(0,g-1,h)&=&e^{-(i+n)x}\left[1+O\left(\frac{nxg}{N}\right)\right]
\\
\Pi(J,g,h)&=&e^{-(i+n)x}\left[1+O\left(\frac{n^2x^2}{g}\right)+
O\left(\frac{nxg}{N}\right)\right].
\end{eqnarray*}
We obtain 
$$
\bigl|{\mathcal B}(i,n) - e^{-(i+n)x}\bigr| \leq \bigl|\Pi(0,g-1,h)-\Pi(J,g,h)\bigr|
$$
and hence
\begin{equation}\label{eq:EstimateB}
{\mathcal B}(i,n) = e^{-(i+n)x}\left[1+O\left(\frac{nxg}{N}\right)
   +O\left(\frac{n^2x^2}{g}\right)\right].
\end{equation}
This estimate is valid for both odd and even period.

Next we consider the products ${\mathcal A}(i,n)$ in (\ref{eq:Ao}) and (\ref{eq:Ae}).
Recalling the definition of $f(N)$ given in equation (\ref{eq:f}), and the fact 
that $h=O(g)$, we have the following asymptotic behaviour
\begin{eqnarray}
{\mathcal A}(0,0)=1
&&\qquad
{\mathcal A}(i,n)=f(N)^{i+n}
    \left[1+O\left(\frac{xn^2}{g}\right)\right],\qquad i+n>0 \label{eq:EstimateA}
\end{eqnarray}

We must consider two cases.

\bigskip

\noindent{\sl Case I: $f(N)\not\to 0$.}\quad From (\ref{eq:EstimateB})
and (\ref{eq:EstimateA}), we obtain, for both
even and odd period
\begin{eqnarray}
\mu(t,i)&=&
\frac{1}{i!}\sum_{n\geq0}\frac{(-1)^n}{n!}e^{-(i+n)y}
\left[
1+
O\left(\frac{n^2x}{g}\right)
\right]
\left[
1+
O\left(\frac{nxg}{N}\right)+
O\left(\frac{n^2x^2}{g}\right)
\right]\nonumber\\
&=&
\frac{(e^{-y})^i}{i!}\sum_{n\geq0}\frac{(-1)^ne^{-ny}}{n!}+
\frac{(e^{-y})^i}{i!}\sum_{n\geq0}\frac{(-1)^n}{n!}
\left[
O\left(\frac{n^2x}{g}\right)
+
O\left(\frac{nxg}{N}\right)
\right]\nonumber\\
&=&
\frac{(e^{-y})^i}{i!}
\left[
e^{\displaystyle-(e^{-y})}
+
O\left(\frac{x}{g}\right)
+
O\left(\frac{xg}{N}\right)\label{eq:muError}
\right]
\end{eqnarray}
where we have used the identities
$$
\sum_{n\geq0}\frac{(-1)^nn^2}{n!}=0=O(1)
\hskip 30pt
\sum_{n\geq0}\frac{(-1)^nn}{n!}=-\frac{1}{e}=O(1).
$$

Now we adjust $t$ so that $y$ is constant. This gives $x=O(\ln(f))$, and
for both odd and even period the error terms in (\ref{eq:muError}) vanish.
We have convergence to a Poisson distribution
\begin{equation}\label{eq:Poisson}
\mu(t,i)\sim e^{\displaystyle -\alpha}\frac{\alpha^i}{i!}\qquad
\alpha= e^{-y}.
\end{equation}
Note that for the above result to hold, we have only required that $f(N)\not\to 0$, 
not the existence of a limit for $f(N)$.

If $f(N)$ tends to a non-zero finite limit $c$, there is no logarithmic term 
in the function $y(N)$ ---cf.~(\ref{eq:yNew})--- and it is possible to consider 
the probability $\mu(t,i)$ without scaling the period $t$. 
For fixed $t$, we have $x=O(g/N)\to 0$ and $y\to-\ln(c)$, and equation 
(\ref{eq:muError}) gives
$$
\mu(t,i)\sim e^{-c}\frac{c^i}{i!},
$$
independent of $t$.

\medskip

\noindent{\sl Case II: $f(N)\to0$.}\quad From (\ref{eq:EstimateA}) 
we have ${\mathcal A}(i,n)\sim \delta_{i+n}$ (Kronecker's delta),
and hence
$$
\mu(t,i)\sim\delta_i.
$$
The proof of theorem B is complete.
\endproof

We remark that theorem A can actually be derived from theorem B
since the distribution $\cR(x)$ follows from first principles as
\begin{equation}
\cR(x)=\lim_{N\to\infty}\cR_N(x)=\lim_{N\to\infty} \; \;
 \sum_{t=1}^{\lfloor xz\rfloor }\;  \frac{t}{N} \sum_{i=1}^{\lfloor
N/t\rfloor} i\,  \mu(t,i),
\end{equation}
equivalently $\langle P_t\rangle$ follows from
\begin{equation} \label{eq:NewP}
\langle P_t\rangle = \frac{t}{N} \sum_{i=1}^{\lfloor N/t\rfloor} i\, 
\mu(t,i).
\end{equation}
Note that the variable $x$ in \eqref{eq:yNew} is, asymptotically, the
variable $x$ in $\cR(x)$.
Taking the result \eqref{eq:Poisson} gives
$$
\sum_{i=1}^{\infty} i\,  \mu(t,i) = \alpha\,  e^{\displaystyle
-\alpha}\; \sum_{i=1}^{\infty}   \frac{\alpha^{(i-1)}}{(i-1)!}
=\alpha= e^{-y}
$$
Since $e^{-y}=f\, e^{-x}$ from \eqref{eq:yNew}, we have
\begin{equation} \label{eq:Pgeneral}
\langle P_t\rangle \sim \frac{t}{N}  f\, e^{-x}.
\end{equation}
If $t=2k-1$ is odd, we substitute for $f$ from \eqref{eq:f} and
$x$ from \eqref{eq:yNew} to see
\begin{equation} \label{eq:NewP2}
\langle P_{2k-1}\rangle \sim
\frac{gh (2k-1)}{N^2}\,  e^{-(k-1) \kappa}.
\end{equation}
Comparing this to \eqref{eq:Compare} and noting there that
$\lambda_\pm^{k-1} \sim\lambda^{k-1} = e^{(k-1) \ln(\lambda)} \sim
e^{-(k-1) \kappa}$
shows that we will recover the result of theorem
\ref{theorem:DistributionOdd}.
A similar exercise on \eqref{eq:Pgeneral} for $t=2k$ using the
appropriate $f$ from
\eqref{eq:f} yields the result of theorem
\ref{theorem:DistributionEven} and hence
the distribution $\cR(x)$ upon summation of odd and even
contributions. We have chosen to give the independent
proof of Theorem A in section \ref{section:RandomInvolutions} as it
is simple and direct and bypasses the calculation of $\mu(t,i)$.

The applicability of Theorem B to algebraic maps 
over finite fields should be considered with care. For fixed $t$, 
the maximum number of $t$-cycles for an algebraic map is plainly finite, 
because the points in a cycle are solutions of a set of polynomial 
equations; therefore, in the large field limit, one cannot possibly 
expect convergence to (\ref{eq:muAsymptotic}).
Under what circumstances is then theorem B a model for such maps?

To articulate this question, let us consider a specific example,
namely the H\'enon map (\ref{eq:Henon}), with $a=1$. 
The set $Fix\,G$ is the line $y=x$, and so over the field $\FF_p$, this 
map has parameters $g=h=p$, hence $\alpha=c=1$ in theorem B. 
Let us consider symmetric cycles, which, due to theorem A, are 
expected to dominate the dynamics. When the fixed sets are 
one-dimensional, the study of symmetric orbits is greatly 
simplified, because the points in these orbits are roots of 
polynomial equations in one indeterminate (as opposed to 
being points on varieties). 
Now, any symmetric odd cycle has a point $(x,x)$ on $Fix\,G$, and 
the values of $x$ for the period $t=5$ are found to be the roots of 
the irreducible polynomial
\begin{eqnarray*}
\Phi_5(x)&=&x^{6}-2\,x^{5}+5\,x^{4}-6\,x^{3}+8\,x^{2}-4\,x+3.
\end{eqnarray*}
Thus the number of symmetric $5$-cycles over the field $\FF_p$ is 
equal to the number of linear factors in the factorization of $\Phi_5$ 
modulo $p$. In particular, the maximum number of 5-cycles is equal 
to 6, the degree of $\Phi_5$.

In general, polynomial factorization over a finite field is a highly 
irregular process, subject only to a slight deterministic constraint
---Stickelberger's theorem--- which determines the parity of the number 
of factors \cite[section 4.8]{Cohen}.
Significantly, this process is ruled by a probabilistic law
---Cebotarev's density theorem---
which determines the probability associated with each factorization
type (number and degree of the factors) in terms of properties
of the Galois group of the polynomial \cite[p.~129]{PohstZassenhaus}.
Specifically, for each factorization with distinct roots, the degree of the factors
defines an additive partition of the degree of the polynomial.
Each partition in turn identifies uniquely a conjugacy class of the 
Galois group, with the terms in the partition corresponding to cycle 
lengths in the cycle decompositions of the permutations in the class.
The probability of a factorization type is then given by the
relative size of the corresponding class in the Galois group.
If the Galois group of $\Phi$ is the symmetric group, the largest 
possible one, then every factorization type does actually occur, 
and with positive density among all primes.

It can be verified that the Galois group of $\Phi_5$ is in fact $S_6$, 
the set of all permutations of $\mbox{deg}(\Phi_5)=6$ objects 
\cite{DavenportSmith}. 
There are therefore 11 possible types of factorizations of $\Phi_5$, 
corresponding to the 11 additive partitions of the integer 6.
We factor $\Phi_5(x)$ modulo a few consecutive primes
\begin{equation*}
\Phi_5(x)\,\equiv\,
\begin{cases}
\left(x+508 \right) 
\left(x+917 \right)
\left(x+1165\right)
\left(x+1486 \right)
\left(x+3168 \right) 
\left(x+5880\right)
&\mbox{(mod 6563)}\\
x^{6}+6567\,x^{5}+5\,x^{4}+6563\,x^{3}+8\,x^{2}+6565\,x+3 
&\mbox{(mod 6569)}\\
\left(x+1265\right)
\left(x+3889\right)
\left({x}^{4}+1415\,{x}^{3}+2999\,{x}^{2}+4200\,x+2683 \right) 
&\mbox{(mod 6571).}
\end{cases}
\end{equation*}
These factorizations correspond to the partitions $1+1+1+1+1+1=6=1+1+4$,
respectively. We see that, over the field $\FF_p$, the H\'enon map has 
six symmetric 5-cycles for $p=6563$ (the smallest prime for which this 
happens), no 5-cycles at all for $p=6569$, and two 5-cycles for $p=6571$. 

The degree of the polynomials $\Phi_t$ grows exponentially with $t$, 
and it is reasonable to expect that, typically, they will be irreducible, 
with large Galois groups. 
For the symmetric group, the Cebotarev's probability $\nu(t,i)$ of occurrence 
of $i$ symmetric cycles of period $t$ is just the probability of having exactly 
$i$ fixed points in the appropriate symmetric group, namely
\begin{equation}\label{eq:Cebotarev}
\nu(t,i)=\frac{1}{i!}\sum_{j=0}^{d-i}(-1)^j\frac{1}{j!}
\hskip40pt 
d=\mbox{deg}\,\Phi_t(x).
\end{equation}
In the large $t$ limit, $\nu(t,i)$ converges to $e^{-1}/i!$, independent 
of $t$. This agrees with the prediction of the random involution model.

It is instructive to compare the data (\ref{eq:Cebotarev}) from 
Cebotarev's theorem with the prediction of theorem B, for
the period $t=5$.

\medskip

\hfil\vtop{\baselineskip 15pt
\halign{\quad $#$\hfil &\qquad $\displaystyle#$\hfil &\qquad #\hfil \quad\cr
i&\mbox{Cebotarev thm} & thm B\cr
\noalign{\vskip 5pt\hrule\vskip 5pt}
0& {53}/{144}= 0.3681\ldots & 0.3679\ldots\cr
1& {11}/{30} = 0.3667\ldots & 0.3679\ldots\cr
2& {3}/{16}  = 0.1875\ldots & 0.1839\ldots\cr
3& {1}/{18}  = 0.0555\ldots & 0.0613\ldots\cr
4& {1}/{48}  = 0.0208\ldots & 0.0153\ldots\cr
5&   0                      & 0.0031\ldots\cr
6& {1}/{720} = 0.0014\ldots & 0.0005\ldots\cr
}}\hfil

\medskip
We see that, even for such a low period, for $i<3$ the agreement 
is already reasonable.

Thus, in the two-dimensional case, some statements on the probability 
of multiple occurrences of symmetric $t$-cycles translate into
statements on the maximality of the Galois groups of $\Phi_t$.
Even though this formulation would allow one to re-cast some 
dynamical questions in algebraic terms, we note that problems 
related to Galois groups of iterated polynomials are often 
quite difficult.


\section{Concluding remarks}\label{section:ConcludingRemarks}

When the reduction over a finite field of a reversible algebraic map
with a single family of reversing symmetries is a permutation, 
numerical experiments suggest the period distribution function 
is given by the period distribution function $\cR(x)$ of \eqref{eq:R}.  
In this paper, we have shown that $\cR(x)$ corresponds to the expected 
period distribution function in a probability space
of pairs of involutions characterized by the cardinalities of their fixed sets.
We find that the same probability space furnishes a Poisson law for the probability of 
cycles with the same period, also in agreement with numerics.

We conclude with some remarks:

\vspace{5pt}
\noindent {\bf (1)}
To prove theorem A, only mild regularity conditions
---the existence of the limits (\ref{eq:ghAsymptotics})
--- were imposed
on the behaviour of the sequences $g(N)$, $h(N)$, specifying
the cardinalities of the fixed sets of the two
involutions, as the size $N=q^n$ of the phase space increases
(with $n$ the dimensionality of the map).
Similarly, theorem B involves checking the behaviour of  the sequence
$f(N)$ in \eqref{eq:f}.
We point out that  when a reversible map of a finite set is obtained
by reducing an algebraic mapping to a finite field, the behaviour of such
sequences is strongly constrained, resulting in algebraic growth.

A simple but significant case comprises maps for which the fixed set
$Fix\,G$ is an affine subspace of some linear space.
An example is given by two-dimensional reversible polynomial
automorphisms
over some field $K$, where the existence of a normal form for the two
involutions ensures that their fixed set is a single point or a line
\cite{BaakeRoberts:05}.
Localization to the finite field $\FF_q$ gives a phase space
with $N=q^2$ points, and since a line in the affine plane $\FF_q^2$
has $q$ points, we have an algebraic growth: $g(N)+h(N)=2q=2\sqrt{N}$.
In dimension $d$, each of $Fix\,G$ and $Fix\,H$ can have integer
dimension ranging from 0 to $d-1$, and if these sets are affine spaces, we
obtain the sequence $g(N)+h(N)=N^r+N^s$, with $r,s\in\QQ$, $0\le
r,s<1$.
The conditions \eqref{eq:ghAsymptotics} are then met as long as one
of $r$ or $s$ is nonzero, i.e. as long as one of $Fix\,G$ and $Fix\,H$ is at
least one-dimensional.
Furthermore, in this situation, we see from \eqref{eq:f} that
$f(N)=N^{r+s-1}$ or
$f(N)=(N^{2r-1}+N^{2s-1})/2$. So either $f=1$ (e.g., when $r=s=1/2
\Leftrightarrow g=h=\sqrt{N}$ )
or $f \to\infty$ algebraically or
$f \to 0$ algebraically.

More generally, the sets $Fix\,G$ and $Fix\,H$ are
algebraic varieties over $\FF_q$ determined by $n$ polynomial
equations,
and we are interested in determining the number of points on these
varieties.
The number $M$ of points on an $r$-dimensional irreducible
algebraic variety in the $n$-dimensional projective space of degree
$d$ over
the field $\FF_q$ satisfies the Hasse-Weil bound
\cite{LangWeil}
\begin{equation} \label{eq:HWbound}
|M-q^r|\leq (d-1)(d-2)q^{r-1/2}+Cq^{r-1}
\end{equation}
where $C=C(n,r,d)$.
For $r=1$ we obtain a square root bound, and we can also replace
$(d-1)(d-2)$ by the genus of the curve.
The result \eqref{eq:HWbound} ensures that, in the large field limit,
there is algebraic growth of the cardinality
of fixed sets. So the discussion of the previous paragraph applies
concerning the validity of the asymptotic
conditions  (\ref{eq:ghAsymptotics}) and the behaviour of $f(N)$.

\vspace{5pt}
\noindent {\bf (2)}
Experiments show that algebraic maps without time-reversal symmetry 
behave quite differently from reversible ones, when represented on a 
finite field \cite{RobertsVivaldi:05}.
In the absence of reversibility, the cycle distribution is found to be 
that of a random permutation, given by the identity function on 
the unit interval. This distribution, obtained via the scaling $z(N)=N$,
has led to a general conjecture \cite[conjecture 2]{RobertsVivaldi:05}.

\vspace{5pt}
\noindent {\bf (3)}
When a general rational reversible map is reduced over finite fields,
we expect that points where the denominators vanish (singularities) will
occur. Restricted to its periodic points, the reduced map is still a permutation.
Numerically, we find the distribution of periods is still governed by $\cR(x)$,
consistent with Theorem A  \cite{Jogia,VialletJogiaRobertsVivaldi}.
But nonperiodic orbits will also appear in the reduction of such a map.
We are presently investigating the extension of the random involution model 
to explain various aspects of the dynamics of the reduction to finite fields 
of general rational reversible maps \cite{VialletJogiaRobertsVivaldi}.
We are also extending our numerical experiments to cover the asymptotic
regime $\FF_q$, $q=p^k$, $k\to\infty$.

\bigskip\noindent
{\sc Acknowledgements:} \/ This work was supported by the Australian 
Research Council, via Discovery Project DP0774473, and the
Australian Academy of Science.
JAGR and FV would like to thank, respectively, the School of 
Mathematical Sciences at Queen Mary, University of London, 
and the School of Mathematics and Statistics at the University 
of New South Wales, Sydney, for their hospitality.


\end{document}